\documentclass{article}
\usepackage{amsfonts}
\usepackage{amsmath}

\newtheorem{theorem}{Theorem}

\newtheorem{lemma}[theorem]{Lemma}

\newtheorem{proposition}[theorem]{Proposition}
\newtheorem{remark}[theorem]{Remark}

\newenvironment{proof}[1][Proof]{\noindent\textbf{#1.} }{\ \rule{0.5em}{0.5em}}
\input{tcilatex}
\begin{document}

\title{On the synectic metric in the tangent bundle of a Riemannian manifold}
\author{Melek ARAS\thanks{%
Department of Mathematics,Faculty of Arts and Sciences, Giresun University,
28049, Turkey e-mail:melekaras25@hotmail.com;melek.aras@giresun.edu.tr }}
\maketitle

\begin{abstract}
The purpose of this paper is to investigate applications the covariant
derivatives of the covector fields and killing vector fields with respect to
the\ synectic lift $^{S}g=^{C}g+^{V}a$ in a the Riemannian manifold \ to its
tangent bundle $T\left( M_{n}\right) $ , where $^{C}g$-complete lift of the
Riemannian metric, $^{V}a$-vertical lift of the symmetric tens\"{o}r field
of type $\left( 0,2\right) $ in $M_{n}$ .

\end{abstract}

\begin{description}
\item[1. Introduction] 
\end{description}

Suppose that there is given the following Riemannian metric

\begin{equation}
^{S}\widetilde{g}_{CB}dx^{C}dx^{B}=a_{ji}dx^{j}dx^{i}+2g_{ji}dx^{j}\delta
y^{i}  \tag{$1$}
\end{equation}%
in tangent bundle in $T\left( M_{n}\right) $ over a Riemannian manifold $%
M_{n}$ with metric $g$, where $a_{ji}$ are components of a symmetric tensor
field of type $\left( 0,2\right) $ in $M_{n}$ and $\delta
y^{h}=dy^{h}+\Gamma _{i}^{h}dx^{i}$, $\Gamma _{i}^{h}=y^{j}$ $\Gamma
_{ji}^{h}$ with respect to the induced coordinates $\left(
x^{h},y^{h}\right) $ in $\pi ^{-1}\left( U\right) \subset T\left(
M_{n}\right) $. We call this metric the synectic metric. The synectic metric 
$^{S}g=^{C}g+^{V}a$ has respectively components\cite{2} and its
contravarient components

\begin{equation}
^{S}g=\left( ^{S}\widetilde{g}_{CB}\right) =\left( 
\begin{array}{cc}
a_{ji}+\partial g_{ji} & g_{ji} \\ 
g_{ji} & 0%
\end{array}%
\right) ,^{S}\widetilde{g}^{CB}=\left( 
\begin{array}{cc}
0 & g^{ji} \\ 
g^{ji} & x^{\overline{s}}\partial _{s}g^{ji}-a_{..}^{ji}%
\end{array}%
\right)   \tag{$2$}
\end{equation}%
where $\partial g_{ji}=x^{\overline{s}}\partial _{s}g_{ji}$ and \ $%
a_{..}^{ti}=g^{jt}a_{js}g^{si}$. 

Components of the Riemannian connection determined by the synectic metric $%
^{S}g$ are \cite{1}

\begin{equation}
\left\{ 
\begin{array}{c}
^{S}\Gamma _{ji}^{k}=\Gamma _{ji}^{k},\text{ \ }^{S}\Gamma _{\overline{j}i}^{%
\overline{k}}=\Gamma _{ji}^{k},\text{ \ \ \ }^{S}\Gamma _{j\overline{i}}^{%
\overline{k}}=\Gamma _{ji}^{k},\text{\ \ }^{S}\Gamma _{\overline{j}\overline{%
i}}^{\overline{k}}=0 \\ 
^{S}\Gamma _{j\overline{i}}^{k}=^{S}\Gamma _{\overline{j}i}^{k}=^{S}\Gamma _{%
\overline{j}\overline{i}}^{k}=0,\text{ \ \ }^{S}\Gamma _{ji}^{\overline{k}%
}=x^{\overline{t}}\partial _{t}\Gamma _{ji}^{k}+H_{ji}^{k}%
\end{array}%
\right.  \tag{$3$}
\end{equation}%
with respect to the induced coordinates in $T\left( M_{n}\right) $, $\Gamma
_{ji}^{k}$ being Christoffel symbols constructed with $\ g_{ji}$ . Where $\
H_{ji}^{k}=\frac{1}{2}g^{ks}\left( \nabla _{j}a_{si}+\nabla
_{i}a_{js}-\nabla _{s}a_{ji}\right) $ is a tensor of type $\left( 1,2\right) 
$ and $\nabla _{s}a_{ji}=\partial _{s}a_{ji}-\Gamma _{kj}^{l}a_{li}-\Gamma
_{ki}^{l}a_{jl}$.

The metric connection $\widetilde{\nabla }$ of the synectic metric satisfies 
$\widetilde{\nabla }_{C}$ $^{S}g_{BA}=0$ and has non-trivial torsion tensor $%
\overline{T}_{CB}^{A}$, which is skew-symmetric in the indices $C$ and $B$
.Then the metric connection $\widetilde{\nabla }$ of the synectic metric \
has components\cite{1}

\begin{equation}
\left\{ 
\begin{array}{l}
\widetilde{\Gamma }_{ji}^{h}=\widetilde{\Gamma }_{\overline{j}i}^{\widetilde{%
h}}=\widetilde{\Gamma }_{j\overline{i}}^{\overline{h}}=\Gamma _{ji}^{h}\text{%
, \ \ \ \ \ \ \ \ \ \ \ \ \ \ \ \ \ \ \ \ \ } \\ 
\widetilde{\Gamma }_{ji}^{h}=\widetilde{\Gamma }_{\overline{j}i}^{\widetilde{%
h}}=\widetilde{\Gamma }_{j\overline{i}}^{\overline{h}}=\Gamma _{ji}^{h}\text{%
, \ \ \ \ \ \ \ \ \ \ \ \ \ \ \ \ \ \ \ \ \ \ } \\ 
\widetilde{\Gamma }_{ji}^{\overline{h}}=x^{\overline{t}}\partial _{t}\Gamma
_{ji}^{h}+H_{ji}^{h}-y^{k}R_{kjih}\text{ \ \ \ \ \ \ \ \ \ \ \ \ \ \ \ \ \ \
\ \ \ }%
\end{array}%
\right.  \tag{$4$}
\end{equation}%
with respect to the induced coordinates, $\Gamma _{ji}^{k}$ being
Christoffel symbols formed with $g_{ji}$, where $H_{ji}^{k}=\frac{1}{2}%
g^{ks}\left( \nabla _{j}a_{si}+\nabla _{i}a_{js}-\nabla _{s}a_{ji}\right) $.

Given a vector field $\widetilde{X}$ in $T\left( M_{n}\right) $, the $1-form$
, $\widetilde{\omega }$ defined by $\widetilde{\omega }\left( \widetilde{Y}%
\right) =\widetilde{g}\left( \widetilde{X},\widetilde{Y}\right) ,\widetilde{Y%
}$ being an arbitrary element of $T_{0}^{1}\left( M_{n}\right) $, is called
the covector field associated with $\widetilde{X}$ and denoted by $%
\widetilde{X}^{\ast }$. If $\widetilde{X}$ has local components $\widetilde{X%
}^{A}$, then the associated covector field $\widetilde{X}^{\ast }$of $%
\widetilde{X}$ has local components $\widetilde{X}_{C}=\widetilde{g}_{CA}%
\widetilde{X}^{A}$.

Let $\omega $be a $1-form$ in $M_{n}$ with components $\omega _{i}$. Then
the vertical , complete and horizontal lifts of $\omega $ to $T\left(
M_{n}\right) $ have respectively components\cite{4}

\begin{equation*}
\left( ^{V}\omega _{B}\right) =\left( \omega _{i},0\right) ,\text{ \ \ \ }%
\left( ^{C}\omega _{B}\right) =\left( \partial \omega _{i},\omega
_{i}\right) ,\text{ \ \ \ }\left( ^{H}\omega _{B}\right) =\left( -\Gamma
_{i}^{k}\omega _{k},\omega _{i}\right)
\end{equation*}%
with respect to the induced coordinates in $T\left( M_{n}\right) .$

The associated covector fields of the vertical , complete and horizontal
lifts to $T\left( M_{n}\right) $, with the synectic metric , of a vector
field $X$ with components $X^{h}$ in $M_{n}$ are respectively

\begin{equation}
\left( X_{i},0\right) ,\text{ \ \ }\left( y^{s}\partial
_{s}X_{i}+a_{ij}X^{j},X_{i}\right) ,\text{ \ }\left( \Gamma
_{i}^{h}X_{h}+a_{ij}X^{j},X_{i}\right)  \tag{$5$}
\end{equation}%
with respect to the induced coordinates, where $X_{j}=g_{ji}X^{i}$ are
components of the covector field $X^{\ast }$ associated with $X$.

A vector field $X\epsilon \Im _{0}^{1}$ $\left( M_{n}\right) $ is said to be
a Killing vector field of a Riemannian manifold with metric $g$, if $%
\tciLaplace _{X}g=0\cite{3}$. In terms of components $g_{ji}$ of $g$, $X$ is
a Killing vector field if and only if

\begin{center}
$\tciLaplace _{X}g=X^{\alpha }\nabla _{\alpha }g_{ji}+g_{\alpha i}\nabla
_{j}X^{\alpha }+g_{j\alpha }\nabla _{i}X^{\alpha }=\nabla _{j}X_{i}+\nabla
_{i}X_{j}=0$,
\end{center}

$X^{\alpha }$ being components of $X$, where $\nabla $ is the Riemannian
connection of the metric $g$.

\begin{description}
\item[2. Main Results] 
\end{description}

We now take a vector field $X$ in $M_{n}$ with components $X^{h}$.Then,
since the associated covector fields of the lifts of $X$ have respectively
components give by $\left( 5\right) $, we have by $\left( 3\right) $ and $%
\left( 5\right) $

\begin{equation}
\left\{ 
\begin{array}{l}
\begin{array}{l}
^{S}\nabla _{B}^{V}X_{A}=\left( 
\begin{array}{cc}
\nabla _{j}X_{i} & 0 \\ 
0 & 0%
\end{array}%
\right) , \\ 
^{S}\nabla _{B}^{C}X_{A}=\left( 
\begin{array}{cc}
\partial \left( \nabla _{j}X_{i}\right) +\nabla _{j}\left(
a_{il}X^{l}\right) -H_{ji}^{m}X_{m} & \nabla _{j}X_{i} \\ 
\nabla _{j}X_{i} & 0%
\end{array}%
\right)%
\end{array}
\\ 
\end{array}%
\right.  \tag{6}
\end{equation}%
and consequently

\begin{equation}
\left\{ 
\begin{array}{l}
^{S}\nabla _{B}^{V}X_{A}+^{S}\nabla _{A}^{V}X_{B}=\left( 
\begin{array}{cc}
\nabla _{j}X_{i}+\nabla _{i}X_{j} & 0 \\ 
0 & 0%
\end{array}%
\right) , \\ 
^{S}\nabla _{B}^{C}X_{A}+^{S}\nabla _{A}^{C}X_{B}=\left( 
\begin{array}{cc}
^{S}\nabla _{j}^{C}X_{i}+^{S}\nabla _{i}^{C}X_{j} & ^{S}\nabla _{j}^{C}X_{%
\overline{i}}+^{S}\nabla _{\overline{i}}^{C}X_{j} \\ 
^{S}\nabla _{\overline{j}}^{C}X_{i}+^{S}\nabla _{i}^{C}X_{\overline{j}} & 
^{S}\nabla _{\overline{j}}^{C}X_{\overline{i}}+^{S}\nabla _{\overline{i}%
}^{C}X_{\overline{j}}%
\end{array}%
\right) \\ 
^{S}\nabla _{j}^{C}X_{i}+^{S}\nabla _{i}^{C}X_{j}=\partial \left( \nabla
_{j}X_{i}+\nabla _{i}X_{j}\right) \\ 
\text{ \ \ \ \ \ \ \ \ \ \ \ \ \ \ \ \ \ }+\nabla _{j}\left(
a_{il}X^{l}\right) +\nabla _{i}\left( a_{jl}X^{l}\right)
-H_{ji}^{m}X_{m}-H_{ij}^{m}X_{m} \\ 
\begin{array}{c}
^{S}\nabla _{j}^{C}X_{\overline{i}}+^{S}\nabla _{\overline{i}%
}^{C}X_{j}=\nabla _{j}X_{i}+\nabla _{i}X_{j} \\ 
^{S}\nabla _{\overline{j}}^{C}X_{i}+^{S}\nabla _{i}^{C}X_{\overline{j}%
}=\nabla _{j}X_{i}+\nabla _{i}X_{j}^{S} \\ 
\nabla _{\overline{j}}^{C}X_{\overline{i}}+^{S}\nabla _{\overline{i}}^{C}X_{%
\overline{j}}=0%
\end{array}
\\ 
\\ 
\end{array}%
\right.  \tag{7}
\end{equation}%
with respect to the induced coordinates ,where $X_{j}=g_{jk}X^{k}$. From $%
\left( 7\right) $ we have

\begin{theorem}
Necessery and sufficient conditions in order that $\left( a\right) $ the
vertical, $\left( b\right) $ complete lifts to $T\left( M_{n}\right) $, with
the synectic metric, of a vector field $X$ in $M_{n}$ be a Killing vector
field \ in $T\left( M_{n}\right) $ are that, respectively, $\left( a\right) $
$X$ is a Killing vector field in $M_{n}$and $\left( b\right) $ $X$ is
Killing vector field with vanishing covariant \ derivative in $M_{n}$ and
the covariant \ derivative of symmetric tensor field a of type $\left(
0,2\right) $ vanishes.
\end{theorem}

We also have by $\left( 2\right) $ and $\left( 6\right) $

\begin{equation}
\left\{ 
\begin{array}{l}
^{S}\nabla _{B}^{V}X_{A}-^{S}\nabla _{A}^{V}X_{B}=\left( 
\begin{array}{cc}
\nabla _{j}X_{i}-\nabla _{i}X_{j} & 0 \\ 
0 & 0%
\end{array}%
\right) , \\ 
^{S}\nabla _{B}^{C}X_{A}-^{S}\nabla _{A}^{C}X_{B}=\left( 
\begin{array}{cc}
^{S}\nabla _{j}^{C}X_{i}-^{S}\nabla _{i}^{C}X_{j} & ^{S}\nabla _{j}^{C}X_{%
\overline{i}}-^{S}\nabla _{\overline{i}}^{C}X_{j} \\ 
^{S}\nabla _{\overline{j}}^{C}X_{i}-^{S}\nabla _{i}^{C}X_{\overline{j}} & 
^{S}\nabla _{\overline{j}}^{C}X_{\overline{i}}-^{S}\nabla _{\overline{i}%
}^{C}X_{\overline{j}}%
\end{array}%
\right) \\ 
^{S}\nabla _{j}^{C}X_{i}-^{S}\nabla _{i}^{C}X_{j}=\partial \left( \nabla
_{j}X_{i}-\nabla _{i}X_{j}\right) \\ 
\text{ \ \ \ \ \ \ \ \ \ \ \ \ \ \ \ \ \ }+\nabla _{j}\left(
a_{il}X^{l}\right) -\nabla _{i}\left( a_{jl}X^{l}\right)
-H_{ji}^{m}X_{m}+H_{ij}^{m}X_{m} \\ 
\begin{array}{c}
^{S}\nabla _{j}^{C}X_{\overline{i}}+^{S}\nabla _{\overline{i}%
}^{C}X_{j}=\nabla _{j}X_{i}-\nabla _{i}X_{j} \\ 
^{S}\nabla _{\overline{j}}^{C}X_{i}+^{S}\nabla _{i}^{C}X_{\overline{j}%
}=\nabla _{j}X_{i}-\nabla _{i}X_{j}^{S} \\ 
\nabla _{\overline{j}}^{C}X_{\overline{i}}-^{S}\nabla _{\overline{i}}^{C}X_{%
\overline{j}}=0%
\end{array}
\\ 
\\ 
\end{array}%
\right.  \tag{8}
\end{equation}%
with respect to the induced coordinates ,where $^{S}g^{BA\text{ }S}\nabla
_{B}^{V}X_{A}=0$ and $^{S}g^{BA\text{ }S}\nabla _{B}^{C}X_{A}=2g^{ji}\nabla
_{j}X_{i}$. Thus we have, from $\left( 8\right) $, respectively,

\begin{theorem}
The vertical lift of a vector field in $M_{n}$ to $T\left( M_{n}\right) $
with the synectic metric $^{S}g$ is harmonic if and only if the vector field
in $M_{n}$ is closed.
\end{theorem}

\begin{theorem}
The complete lift of a vector field in $M_{n}$ to $T\left( M_{n}\right) $
with the synectic metric $^{S}g$ is harmonic if and only if the vector field
in $M_{n}$ is harmonic and the covariant \ derivative of symmetric tensor
field a of type $\left( 0,2\right) $ vanishes.
\end{theorem}

We consider a vector field $X\epsilon \Im _{0}^{1}\left( M_{n}\right) $.
Then its vertical , complete and horizontal lifts have components of the form

\begin{equation}
^{V}X=\left( 
\begin{array}{c}
0 \\ 
X^{h}%
\end{array}%
\right) ,\text{ \ \ }^{C}X=\left( 
\begin{array}{c}
X^{h} \\ 
\partial X^{h}%
\end{array}%
\right) \text{, \ \ }^{H}X=\left( 
\begin{array}{c}
X^{h} \\ 
-\Gamma _{i}^{h}X^{i}%
\end{array}%
\right)  \tag{9}
\end{equation}%
with respect to the induced coordinates in $T\left( M_{n}\right) $, where $%
\Gamma _{i}^{h}X^{i}=y^{s}$ $\Gamma _{si}^{h}X^{i}.$

Let $X$ be a vector field in $M_{n}$ with local components $X^{k}$. Then,
from $\left( 9\right) $ and $\left( 4\right) $, we see that, the covariant
derivatives of the vertical, complete and horizontal lifts of $X\epsilon \Im
_{0}^{1}\left( M_{n}\right) $with the metric connection $\widetilde{\nabla }$
have respectively components

\begin{equation}
\left\{ 
\begin{array}{c}
\widetilde{\nabla }_{B}\text{ }^{V}X^{A}=\left( 
\begin{array}{cc}
0 & 0 \\ 
\nabla _{j}X^{h} & 0%
\end{array}%
\right) , \\ 
\widetilde{\nabla }_{B}\text{ }^{C}X^{A}=\left( 
\begin{array}{cc}
\nabla _{j}X^{h} & 0 \\ 
\partial \left( \nabla _{j}X^{h}\right)
+H_{jm}^{h}X^{m}-y^{k}R_{kjm}^{h}X^{m} & \nabla _{j}X^{h}%
\end{array}%
\right) , \\ 
\widetilde{\nabla }_{B}\text{ }^{H}X^{A}=\left( 
\begin{array}{cc}
\nabla _{j}X^{h} & 0 \\ 
-\Gamma _{i}^{h}\left( \nabla _{j}X^{i}\right) +H_{jm}^{h}X^{m} & 0%
\end{array}%
\right)%
\end{array}%
\right.  \tag{10}
\end{equation}%
with respect to the induced coordinates $T\left( M_{n}\right) .$

\begin{remark}
$\widetilde{\nabla }=\overline{\nabla }+^{V}H$, where $\overline{\nabla }$
is the metric connection with the metric $^{C}g\cite{1}$
\end{remark}

\begin{remark}
The metric connection $\overline{\nabla }$ coincides with the horizontal
lift $^{H}\nabla $ of Levi Civita connection $\nabla $of $g$ in $M_{n}\cite%
{4}$ . Thus we have
\end{remark}

\begin{proposition}
Necessery and sufficient conditions in order that $\left( a\right) $ the
vertical, $\left( b\right) $ complete and horizontal lifts of a vector field
in $M_{n}$ to $T\left( M_{n}\right) $ with the metric connection $\widetilde{%
\nabla }$ be parallel in $T\left( M_{n}\right) $ are that, respectively, $%
\left( a\right) $ the vector field given in $M_{n}$ is parallel $\left(
b\right) $ the vector field given in $M_{n}$ is parallel and the covariant \
derivative of symmetric tensor field a of type $\left( 0,2\right) $ vanishes.
\end{proposition}

Since $\nabla _{j}X^{h}=t\delta _{j}^{h}$ with constant $t$ implies $%
R_{kji}^{h}X^{i}=0$, we have also

\begin{proposition}
The complete lift of a vector in $M_{n}$ to $T\left( M_{n}\right) $ with the
metric connection $\widetilde{\nabla }$is concurrent if and only if the
vector field given in $M_{n}$ is concurrent and the covariant \ derivative
of symmetric tensor field a of type $\left( 0,2\right) $ vanishes.
\end{proposition}

A vector field $X\epsilon \Im _{0}^{1}$ $\left( M_{n}\right) $ is said to be
an infinitesimal isometry or a Killing vector field of a Riemannian manifold
with metric $g$, if $\tciLaplace _{X}g=0$. In terms of components $g_{ji}$
of $g$, $X$ is an infinitesimal isometry if and only if

\begin{equation}
X^{\gamma }\partial _{\gamma }g_{\alpha \beta }+g_{\alpha \gamma }\partial
_{\beta }X^{\gamma }+g_{\gamma \beta }\partial _{\alpha }X^{\gamma }=0, 
\tag{11}
\end{equation}

$X^{\alpha }$ being components of $X\cite{3}$, where the indices $\alpha
,\beta $ and $\gamma $ run over the range $\left\{ 1,2,...,m\right\} .$

Let there be given in $M_{n}$ a Riemannian metric $g$ with components $%
g_{ji}.$Let $\widetilde{X}$ be a vector field with components $\left( 
\begin{array}{c}
\widetilde{X}^{k} \\ 
\widetilde{X}^{\overline{k}}%
\end{array}%
\right) $ with respect to the induced coordinates in $T\left( M_{n}\right) .$

with respect to the induced coordinates in $T\left( M_{n}\right) $. Then,
taking account of $\left( 2\right) $, we see by virtue of $\left( 11\right) $
that $\widetilde{X}$ is an  infinitesimal isometry n $T\left( M_{n}\right) $
with metric $^{S}g$ if and only if 

\begin{equation}
\left\{ 
\begin{array}{l}
\left( \widetilde{X}^{h}\partial _{h}\partial g_{ji}+\widetilde{X}^{%
\overline{h}}\partial _{h}g_{ji}\right) +\left( \partial g_{jh}\partial _{i}%
\widetilde{X}^{h}+g_{jh}\partial _{i}\widetilde{X}^{\overline{h}}\right)  \\ 
+\left( \partial g_{hi}\partial _{j}\widetilde{X}^{h}+g_{hi}\partial _{j}%
\widetilde{X}^{\overline{h}}\right) +\left( \widetilde{X}^{h}\partial
_{h}a_{ji}+a_{jh}\partial _{i}\widetilde{X}^{h}+a_{hi}\partial _{j}%
\widetilde{X}^{h}\right) =0%
\end{array}%
\right.   \tag{12}
\end{equation}

\begin{equation}
\widetilde{X}^{h}\partial _{h}g_{ji}+g_{jh}\partial _{i}\widetilde{X}%
^{h}+\left( \partial g_{hi}\partial _{\overline{j}}\widetilde{X}%
^{h}+g_{hi}\partial _{\overline{j}}\widetilde{X}^{\overline{h}}\right)
+a_{hi}\partial _{\overline{j}}\widetilde{X}^{h}=0  \tag{13}
\end{equation}

\begin{equation}
\widetilde{X}^{h}\partial _{h}g_{ji}+g_{hi}\partial _{j}\widetilde{X}%
^{h}+\left( \partial g_{jh}\partial _{\overline{i}}\widetilde{X}%
^{h}+g_{jh}\partial _{\overline{i}}\widetilde{X}^{\overline{h}}\right)
+a_{jh}\partial _{\overline{i}}\widetilde{X}^{h}=0  \tag{14}
\end{equation}

\begin{equation}
g_{jh}\partial _{\overline{i}}\widetilde{X}^{h}+g_{hi}\partial _{\overline{j}%
}\widetilde{X}^{h}=0.  \tag{15}
\end{equation}

We shall now prove

\begin{lemma}
Let $C$ be an element of $\Im _{1}^{1}\left( M_{n}\right) $. Then  $%
\tciLaplace _{\iota C}$ $^{S}g=0$ holds if and only if $C=0.$

\begin{proof}
Denote by $C_{i}^{k}$ the local components of $C$. Then $\iota C$ has
components $\left( 
\begin{array}{c}
0 \\ 
y^{i}C_{i}^{k}%
\end{array}%
\right) $ with respect to the induced coordinates in $T\left( M_{n}\right) .$%
Thus, substituting $\widetilde{X}^{k}=0$ and $\widetilde{X}^{\overline{k}%
}=y^{i}C_{i}^{k}$ in $\left( 13\right) $ , we have $g_{ki}C_{j}^{k}=0$,
which implies $C_{j}^{k}=0$, i.e., $C=0$. 
\end{proof}
\end{lemma}

\end{document}